\documentclass[a4paper,titlepage,oneside]{article}

\usepackage{arxiv}

\usepackage[utf8]{inputenc} 
\usepackage[T1]{fontenc}    
\usepackage{hyperref}       
\usepackage{url}            
\usepackage{booktabs}       
\usepackage{amsfonts}       
\usepackage{nicefrac}       
\usepackage{microtype}      
\usepackage{lipsum}		
\usepackage{graphicx}
\usepackage{natbib}
\usepackage{doi}
\usepackage{graphicx}
\usepackage{mathptmx}
\usepackage{amsmath}
\usepackage{amsthm}
\usepackage{amssymb}
\usepackage{multirow}
\usepackage{mathtools}
\usepackage{siunitx}
\usepackage{float}
\usepackage{slashbox}
\usepackage[ruled,vlined,noline]{algorithm2e}
\newtheorem{theorem}{Theorem}

\usepackage{perpage} 
\MakePerPage{footnote} 

\DeclareMathOperator{\prox}{prox}

\newcommand{\R}{\mathbb{R}}
\newcommand{\N}{\mathbb{N}}
\renewcommand{\div}{\operatorname{div}}
\newcommand{\D}{D}
\newcommand{\Dtwo}{D^2}
\newcommand{\Div}{\div}
\newcommand{\Divtwo}{\div^2}
\newcommand{\domain}{\mathcal{M}}

\newcommand{\Kop}{K}
\newcommand{\Lop}{L}
\newcommand{\Epsop}{E}

\renewcommand{\d}[1]{\ensuremath{\operatorname{d}\!{#1}}}
\newcommand{\bbar}[1]{#1}

\usepackage[acronym]{glossaries}
\newacronym{tgv}{TGV}{total generalized variation}
\newacronym{tv}{TV}{total variation}
\newacronym{mse}{MSE}{mean squared error}
\newacronym{psnr}{PSNR}{peak signal-to-noise ratio}
\newacronym{ssim}{SSIM}{structural similarity index measure}

\title{Learned Discretization Schemes for the Second-Order Total Generalized Variation}

\author{{Lea Bogensperger}\\
	\texttt{lea.bogensperger@icg.tugraz.at} \\
	\And
	{Antonin Chambolle} \\
\texttt{antonin.chambolle@ceremade.dauphine.fr} \\
\And 
{Alexander Effland} \\
\texttt{effland@iam.uni-bonn.de} \\
\And 
{Thomas Pock} \\
\texttt{pock@icg.tugraz.at}
}

\renewcommand{\shorttitle}{Learned Discretization Schemes for the Second-Order Total Generalized Variation}

\hypersetup{
pdftitle={Learned Discretization Schemes for the Second-Order Total Generalized Variation},
pdfauthor={Lea Bogensperger, Antonin Chambolle, Alexander Effland, Thomas Pock},
pdfkeywords={Total generalized variation, discretization, image denoising, bilevel optimization, piggyback algorithm, learning, primal-dual algorithm},
}

\begin{document}
\maketitle

\begin{abstract}
	The total generalized variation extends the total variation by incorporating higher-order smoothness.
Thus, it can also suffer from similar discretization issues related to isotropy. 
Inspired by the success of novel discretization schemes of the total variation, there has been recent work to improve the second-order total generalized variation discretization, based on the same design idea. 
In this work, we propose to extend this to a general discretization scheme based on interpolation filters, for which we prove variational consistency.
We then describe how to learn these interpolation filters to optimize the discretization for various imaging applications.
We illustrate the performance of the method on a synthetic data set as well as for natural image denoising.
\end{abstract}


\section{Introduction}
The \gls{tv} is a popular regularizer for many tasks in image reconstruction, yet it assumes as a prior that images/signals are essentially piecewise constant. 
The extension known as \gls{tgv}~\cite{BrKu10} is a natural way to incorporate more complex signals (such as affine) in the prior, by combining the \gls{tv} of different orders of derivatives into a global image.
Like \gls{tv}, \gls{tgv} can be used in a plug-and-play style in various inverse problems~\cite{knoll2011tgv,niu2014sparse,ranftl2013minimizing,huber2019total}. 
In the continuous domain, it reads as 
\begin{equation}\label{eq:tgv_cont_dual}
    \text{TGV}_{\alpha}^2(u)= \sup_{p} \bigg \{ \int_\Omega u ~\div^2 p \d{x}: p \in \mathcal{C}^\infty(\Omega, \text{Sym}^{2\times 2}), \Vert p \Vert_{\infty}\leq \alpha_0, \Vert \div p\Vert_{\infty} \leq \alpha_1 \bigg \},
\end{equation}
where $\text{Sym}^{2\times 2}$ denotes the space of second-order symmetric tensors, $\alpha=(\alpha_0, \alpha_1)$ are positive parameters, $\div$ denotes the row-wise (or column-wise, as $p$ is symmetric) divergence, and $\div^2$ the divergence of the resulting vector. 
Note that one can similarly define regularizers combining higher orders of derivatives, yet the most common version used is $\text{TGV}_\alpha^2$, therefore we stick to this case. 
Like \gls{tv}, \gls{tgv} is
difficult to discretize while preserving isotropy and rotational invariance.  
For \gls{tv}, improved discretization schemes have been studied in earlier works such as~\cite{condat2017discrete,chambolle2021learning}. 
Recently, a discretization scheme was proposed in~\cite{hosseini2022second} to improve second-order \gls{tgv} inspired by the work of Condat~\cite{condat2017discrete}.
The idea is to impose the constraints in~\eqref{eq:tgv_cont_dual} on the dual variables in the discretized setting on $n$ times denser grids using interpolations to deal with staggered pixel grids arising from the discretized finite difference operators. 

We would like to extend on this work by expressing it in a more general framework for which we show consistency (see Theorem~\ref{th:consistency}). 
This framework is based on local interpolation operations and requires only bounded filter kernels.
Moreover, since it is not straightforward how to choose ideal filters, and this
may depend on the underlying data and the context of the inverse problem, the question arises whether this can be further improved.
An appealing idea is therefore to resort to learning such interpolation filters and subsequently investigate their performance, as recently done in~\cite{chambolle2021learning} for \gls{tv}. 

\section{Problem Setting}
\subsection{Notation}
Let $M,N$ be the dimension of the pixel grid. We usually denote an image $u\in \domain \coloneqq \R^{M\times N}$. For convenience, but with a slight abuse of notation, note that $\domain$ determines the size of a pixel grid, whilst not assuming anything on the respective spatial locations within the grid. 

If no specific norm is indicated, the $\Vert x \Vert_{1,2}$ norm is assumed, which is for $x \in \domain^J$ the absolute sum of the 2-norm of its $J$ components.
Further, we set $\Vert x \Vert_Z=\Vert x \Vert_{1,1,2}$ for $x \in \domain^{J\times I}$, which is the absolute sum consisting of $I$ components of the 2-norm of its $J$ components. 
Finally, let $\Vert \cdot \Vert_Z^\ast$ denote its corresponding dual norm. 

\subsection{Finite Difference Operators}
On a standard Euclidean grid of size $hM\times hN$ with $M\times N$ pixels of size $h\times h$, we define the discrete forward operator $\D:\domain \to \domain^{2}$ for $u\in\domain$ via
$D u = ((Du)^1, (Du)^2) $, where
\begin{align*}
    &(Du)^1_{i+\frac12,j} = \tfrac{1}{h}(u_{i+1,j} - u_{i,j}) \quad &i \leq M-1, j\leq N, \\
    &(Du)^2_{i,j+\frac12} = \tfrac{1}{h}(u_{i,j+1} - u_{i,j})\quad &i \leq M, j\leq N-1.
\end{align*}
To ease the notation we set the values of the derivatives to 0 using Neumann boundary conditions if the index dies out before reaching $M$ or $N$, which also implies that our resulting pixel grids remain of the same size. The tensor-valued symmetric counterpart is given by $\Epsop: \domain^2 \to \domain^3$, and its individual operator components also consist of forward differences\footnote{Note that one could also resort to backward differences, however, for designing suitable interpolation operators on the dual variables this scheme is more convenient since it leads to the component consisting of mixed derivatives being located at the pixel corner.}. Therefore, for $w=(w^1_{i+\frac12,j},w^2_{i,j+\frac12})$ one obtains the symmetrized tensor field
$\Epsop w = \begin{pmatrix}
    (\Epsop w)^1 &  (\Epsop w)^2 \\ (\Epsop w)^2 &  (\Epsop w)^3
\end{pmatrix}$
with 
\begin{align*}
&(\Epsop w)^1_{i+1,j} = \tfrac1h (w^1_{i+\frac32,j}-w^1_{i+\frac12,j})\qquad & i \leq M-1, j\leq N, \\
&(\Epsop w)^2_{i+\frac12,j+\frac12} = \tfrac{1}{2h} (w^1_{i+\frac12,j+1}-w^1_{i+\frac12,j}+w^2_{i+1,j+\frac12}-w^2_{i,j+\frac12}) \qquad &i \leq M-1, j\leq N-1,\\
&(\Epsop w)^3_{i,j+1}  = \tfrac1h (w^2_{i,j+\frac32} -w^2_{i,j+\frac12})\qquad &i \leq M, j\leq N-1.
\end{align*}
Again, the same handling of derivatives using Neumann boundary conditions is used. The symmetrized second-order finite difference operator is then given by $\Dtwo=\Epsop\D$. The corresponding adjoint operators $\Div$ and $\Divtwo$ are directly given by the discrete Gauss-Green theorem. 

\subsection{Second-Order \gls{tgv} Discretization}
\begin{figure}
    \centering
    \includegraphics[width=\textwidth]{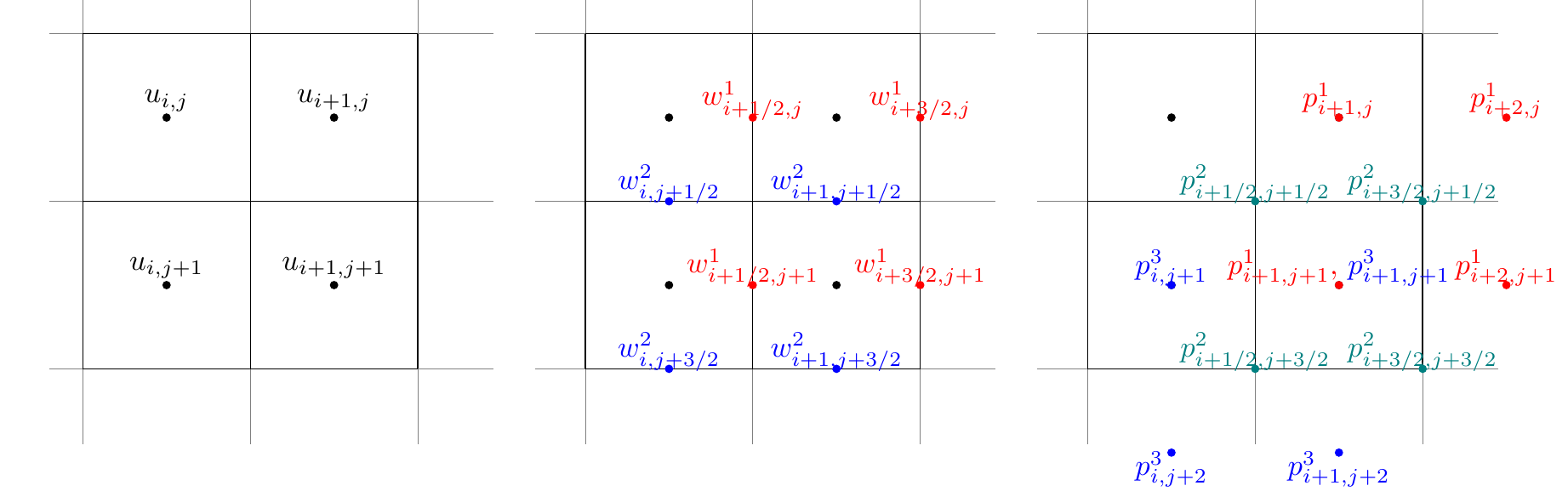}
    \caption{Resulting pixel grids for $w$ and $p$ given an input $u$. Colors indicate different components of the vector/tensor fields for visualization purposes (best viewed on screen).}
    \label{fig:pixelgrid}
\end{figure}

In the spirit of the recently proposed discretization~\cite{hosseini2022second} that builds upon the ideas of Condat's discretization~\cite{condat2017discrete}, the aim is to state a more generalized definition of second-order \gls{tgv} using interpolation filters $\Kop$ and $\Lop$
\begin{align} \label{eq:filters_def}
\Kop  \colon \domain^3 \to \domain^{3 \times n_K}, \text{ } \Lop   \colon  \domain^{2} \to \domain^{2\times n_L},
\end{align} 
with $n_K$ and $n_L$ denoting the number of filters, respectively. These filters can be chosen according to~\cite{hosseini2022second} as shown in Figure~\ref{fig:handcrafted}, but also other choices exist (such as interpolating to arbitrary pixel grid locations), all being based on a staggered grid discretization. 
We start from the standard second-order \gls{tgv} discretization in the primal domain:
\begin{equation}
    \min_{w \in \domain^2} \alpha_{1} \Vert \D u-w\Vert + \alpha_0 \Vert \Epsop w \Vert.
\end{equation}
Using interpolation filters from~\eqref{eq:filters_def} this can be rewritten with $v_K \in \domain^{3 \times n_K}$ and $v_L \in \domain^{2\times n_L}$ as
\begin{align}
     \min_{v_K, v_L, w} \alpha_1 \Vert v_L\Vert_Z + \alpha_0 \Vert v_K \Vert_Z, \text{ s.t. }
    \D u-w = \Lop^\ast v_L, \text{ }     \Epsop w = \Kop^\ast v_K,
\end{align}
where $w$ can be eliminated from the constraints such that we obtain 
\begin{align}\label{eq:tgv_disc_primal}
    \min_{v_K, v_L} \alpha_1 \Vert v_L\Vert_Z + \alpha_0 \Vert v_K \Vert_Z, \text{ s.t. }
    \Dtwo u = \Epsop \Lop^\ast v_L + \Kop^\ast v_K. 
\end{align}
In this sense, one possible interpretation is that we seek to learn a group-sparse coding for the symmetrized second-order discrete derivatives $\Dtwo$ of $u$.
For smooth regions, $v_L$ will be close to 0 and the second-order gradients of the image will only be given by $\Kop^\ast v_K$, whereas for discontinuities $v_L$ contributes as well.
Since $E \Lop^\ast v_L$ will be symmetrized, it essentially leaves more freedom to the model as only the symmetric part of the second-order derivatives in the constraint must be fulfilled. 
Using convex conjugates and duality, the corresponding dual problem reads as
\begin{align}
    \sup_{p\in \domain^3} \langle \Dtwo u, p\rangle, \text{ s.t. }
    \Vert \Lop \Div p \Vert_Z^\ast  \leq \alpha_1, \text{ }
    \Vert \Kop p \Vert_Z^\ast \leq \alpha_0.
\end{align}
Figure~\ref{fig:pixelgrid} shows the resulting pixel grids for the vector and tensor fields $w$ and $p$ arising from the finite difference operators $\D$, $E$, and $\Dtwo$. This basically suggests considering four different pixel locations for interpolation: the pixel center, the center of the horizontal and vertical edges, and the corner.
All other pixel positions at this scale are contained in a superset of these four positions.

\subsection{Interpolation Operators}
\begin{figure}
    \centering
    \includegraphics[width=\textwidth]{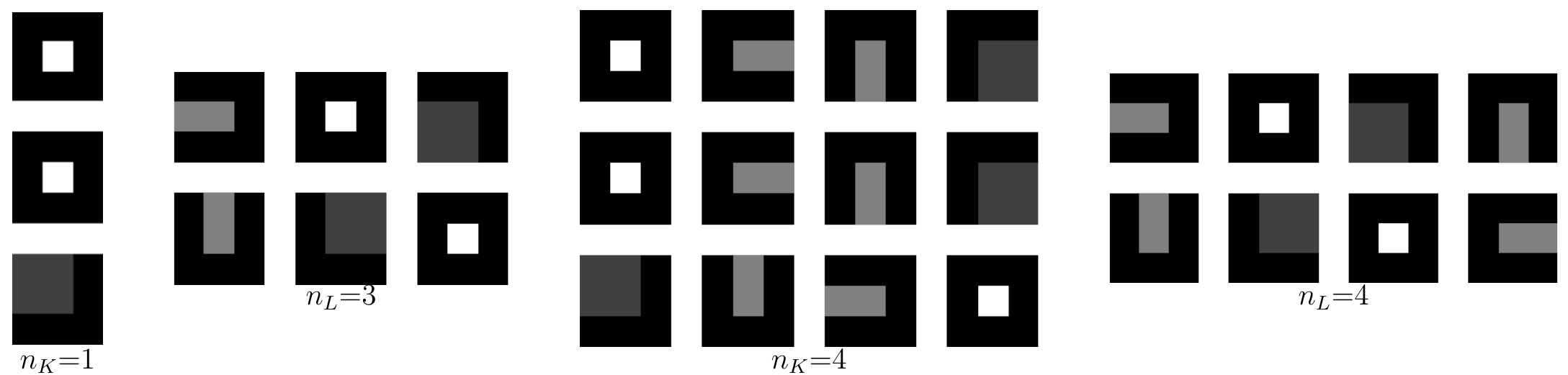}
    \caption{Handcrafted interpolation filters where the intensity values are in $\{0,0.25,0.5,1\}$ (ranging from dark to light) to ensure normalized filter coefficients. The authors in~\cite{hosseini2022second} use $n_K$=1 for $K$ and $n_L$=3 for $L$ ensuring that each component is interpolated from its resulting pixel grid location given in Figure~\ref{fig:pixelgrid} to the pixel center and the horizontal and vertical edges, respectively. In case $n_K$=4 and $n_L$=4 this additionally interpolates to the pixel corner.}
    \label{fig:handcrafted}
\end{figure}

Inspired by the improved discretization schemes on \gls{tv}~\cite{condat2017discrete}, the authors in~\cite{hosseini2022second} construct filters using $n_L=3$ for the dual $\div p$ that is located at the same pixel grid positions as $w$ for both vector field components. Thus both $w^1$, $w^2$ are interpolated to the three pixel grid positions $(i,j),(i+\tfrac12,j),(i,j+\tfrac12)$. While the corresponding interpolation operations are given in detail by~\cite{hosseini2022second}, a schematic of this is also shown in Figure~\ref{fig:handcrafted}.

Thus, the interpolation filter $\big (L^{1,l},L^{2,l} \big )_{l=1}^{n_L}$ is applied using $L^{1,l} w^1 = (L^{1,1} w^1, L^{1,2} w^1, L^{1,3} w^1)$, and analogously for $w^2$. Naturally, this can be extended to also include the fourth position in the pixel corner $(i+\tfrac{1}{2},j+\tfrac{1}{2})$ for $n_L=4$. Using convolutions with filter kernels $(\eta_{m,n}^{1,l},\eta_{m,n}^{2,l})$ and $(\xi_{m,n}^{1,r},\xi_{m,n}^{2,r},\xi_{m,n}^{3,r})$ of limited local support $(2\nu +1)\times(2\nu + 1)$ for $\nu \in \mathbb{N}$ with bounded coefficients, this can be framed in the context of general interpolation operations. The filters can then be expressed as
\begin{align}
    \big (L^l w \big )_{i,j} = \begin{pmatrix} \big (L^{1,l} w^1\big )_{i,j} \\(L^{2,l} w^2\big )_{i,j} 
    \end{pmatrix} =  \begin{pmatrix} \sum_{m,n=-\nu}^\nu \eta_{m,n}^{1,l} w^1_{i+ \frac12 - m,j-n} \\
    \sum_{m,n=-\nu}^\nu \eta_{m,n}^{2,l} w^2_{i - m,j+ \frac12-n}\end{pmatrix}.
\end{align}
In a similar manner, the dual $p$ is interpolated to the pixel position $(i,j)$ for $n_K=1$ for each tensor field component~\cite{hosseini2022second}. Again, the other three positions at both pixel faces and at the corner can be included using $n_K=4$. In the general setting using $\big (K^{1,r},K^{2,r},K^{3,r} \big )_{r=1}^{n_K}$ this amounts to 
\begin{align}
    \big (K^r p \big )_{i,j} = \begin{pmatrix} \big (K^{1,r} p^1\big )_{i,j} \\ (K^{2,r} p^2\big )_{i,j} \\(K^{3,r} p^3\big )_{i,j} 
    \end{pmatrix} =  \begin{pmatrix} \sum_{m,n=-\nu}^\nu \xi_{m,n}^{1,r} p^1_{i- m,j-n} \\
    \sum_{m,n=-\nu}^\nu \xi_{m,n}^{2,r} p^2_{i+\frac12  - m,j+\frac12-n}\\
    \sum_{m,n=-\nu}^\nu \xi_{m,n}^{3,r} p^3_{i - m,j-n}\end{pmatrix}.
\end{align}
In general, however, it is not straightforward how to select the interpolation points within the pixel grid of the dual variables with regards to an improved discretization. 
To gain a basic intuition, experiments on image denoising (see Section~\ref{sec:datasets} on the respective data set and Section~\ref{sec:imagerec}/Algorithm~\ref{alg:pd_general} on the optimization problem/reconstruction algorithm) have been conducted with varying $n_K$ and $n_L$. In general, it seems that the choice of $n_K$ does not impact the performance to a large extent, presumably due to the second-order finite differences which yield very smooth tensor fields. On the other hand, a larger $n_L$ seems to be beneficial, resulting in a denser grid. Moreover, it is noteworthy that these tendencies exhibit small fluctuations depending on the parameters $\alpha_1$ and $\alpha_0$, the choice of data and the level of noise corruption. Due to this ambiguity of selecting the best set of suitable filters, it is tempting to directly learn the filters with the aim to obtain an even better discretization. This is also motivated by the success of learned discretization schemes for \gls{tv}~\cite{chambolle2021learning}.

\section{$\mathrm{\Gamma}$-Convergence of the Discretization}
For simplicity, we use a square grid of $N\times N$ pixels for the domain $\Omega=(0,1)^2$, where each pixel is of size $h\times h$, with $h=1/N$. 
The operators and variables in the discrete setting are now marked with an $h$. We use both primal and dual definitions of the discretized second-order \gls{tgv}
\begin{align}
    \text{TGV}_{\alpha,h}^2(u^h) &= \min_{w^h,v_K^h,v_L^h} \bigg\{ h^2 \alpha_1 \Vert v_L^h\Vert_{Z} + h^2\alpha_0 \Vert v_K^h \Vert_{Z}: \Lop_h^\ast v_L^h = \D_h u^h - w^h \text{, } \Kop_h^\ast v_K^h = \Epsop_h w^h \bigg \} \notag\\&=  \sup_{p^h}\bigg\{h^2 \langle \Div_h^2 p^h, u^h \rangle : \Vert \Lop_h \Div_h p^h \Vert_Z^\ast \leq \alpha_1, \Vert \Kop_h p^h \Vert_Z^\ast \leq \alpha_0 \bigg \}. \label{eq:tgv_disc}
\end{align}
In a slightly simpler setting where we assume that $u$ is global affine plus periodic and $w$ is periodic with periodic boundary conditions, the following theorem states the $\Gamma$-convergence of the discretized second-order \gls{tgv} in~\eqref{eq:tgv_disc} to the continuous second-order \gls{tgv} in~\eqref{eq:tgv_cont_dual}. The corresponding proof is given in the appendix. Note that the minimum in~\eqref{eq:tgv_disc} is attained due to the finite-dimensional setting and the boundedness of $(w,v_K,v_L)$.
\begin{theorem}\label{th:consistency}
We consider the setting
where $u$ is affine plus periodic with period 1 in $\R^2$, and $w$ is $1$-periodic. Then, for interpolation operators $\Kop$ and $\Lop$ that have local support and bounded filter coefficients, 
$\textup{TGV}^2_{\alpha,h}(u^h)$ $\Gamma$-converges to $\textup{TGV}^2_{\alpha}(u)$.
\end{theorem}
The interpretation of this theorem is that
minimizers of problems involving $\text{TGV}^2_{\alpha,h}$ plus some
continuous term (for instance, a quadratic
penalization) will converge, when viewed as piecewise constant functions in the continuum, to minimizers
of the corresponding continuous problem
involving $\text{TGV}^2_{\alpha}$ defined
in~\eqref{eq:tgv_cont_dual}.

\section{Numerical Methods}
\subsection{Image Reconstruction} \label{sec:imagerec}
Second-order \gls{tgv} regularization is typically applied to image reconstruction problems being combined with a task-dependent convex data fidelity term $G(u,f)$, e.g. a typical use case is image denoising with $G(u,f)=\tfrac 1 2 \Vert u-f\Vert_2^2$. Thus, given a corrupted image $f\in \domain$ we obtain the following saddle point problem using the proposed discretization scheme from~\eqref{eq:tgv_disc_primal}
\begin{align}
\min_{u,v_K, v_L} \max_p G(u,f) +\alpha_0 \Vert v_K \Vert_{Z} + \alpha_1 \Vert v_L \Vert_{Z} + \langle \Dtwo u - \Epsop \Lop^\ast v_L - \Kop^\ast v_K, p \rangle.
\label{eq:sp_general}
\end{align} 
This can be solved with a primal-dual algorithm~\cite{chambolle2011first} as described in Algorithm~\ref{alg:pd_general}, using diagonal block-preconditioning~\cite{pock2011diagonal} to determine the step sizes. For details on how to compute the proximal maps, see~\cite{chambolle2021learning}. 
\begin{algorithm}[ht]
\caption{Primal-dual algorithm to solve \eqref{eq:sp_general}.}\label{alg:pd_general}
\SetAlgoLined 
\KwIn{initial values $u^1$, $v_K^1$, $v_L^1$ and $p^1$, block-preconditioned step size parameters  $\tau_u, \tau_{v_K},\tau_{v_L},\sigma>0$, $\theta\in[0,1]$, number of iterations $J$}
\KwResult{approximate saddle point $(u^J, v_K^J, v_L^J, p^J)$}
\For{$j = 1, 2,\ \ldots ,\  J $}
{
$p^{j+1}=p^j + \sigma (\Dtwo u^j - \Epsop \Lop^\ast v_L^j - \Kop^\ast v_K^j)$\;
$\bar{p}^{j+1} = p^{j+1} + \theta( p^{j+1} - p^j)$\;
$u^{j+1} = \prox_{\tau_u G(\cdot,f)}(u^j - \tau_u D^{2^\ast} \bar{p}^{j+1})$\;
$v_L^{j+1} = \prox_{\tau_{v_L} \Vert \cdot \Vert_{Z}} (v_L^j + \tau_{v_L} \Lop \Epsop^\ast \bar{p}^{j+1})$\;
$v_K^{j+1} = \prox_{\tau_{v_K} \Vert \cdot \Vert_{Z}} (v_K^j + \tau_{v_K} \Kop \bar{p}^{j+1})$\;
}
\end{algorithm} 

\subsection{Learning Interpolation Filters} \label{sec:learning}
The interpolation filters can be learned with a bilevel approach, where the outer optimization problem enforces the similarity of the approximate reconstructions $u^\ast$ from the inner problem to a known target data set $t$. 
To achieve this, a loss function is required (we use a quadratic loss $\ell(u^\ast,t)=\tfrac 1 2 \Vert u^\ast -t\Vert_2^2$) with additional constraints on the learned interpolation filters 
\begin{equation}
    \min_{\Kop,\Lop} \frac 1 S \sum_{s=1}^S \ell(u^{s \ast}(\Kop,\Lop),t^s) + \mathcal{R}(\Kop)  + \mathcal{R}(\Lop).
\end{equation}
The constraints on the filters are given by $\mathcal{R}(K)=\delta_{(C_{\Sigma=1})^{3,n_K}}$ and $\mathcal{R}(L)=\delta_{(C_{\Sigma=1})^{2,n_L}}$, with $\delta_{C_{\Sigma=1}}$ the indicator function of the set $C_{\Sigma=1}$ per filter for each component, to ensure the boundedness of the filters such that for each the sum of the coefficients is 1 (also see Section~\ref{sec:filtersettings} for more details). 

As an alternative to an unrolling scheme, we resort to a piggyback-style algorithm for obtaining derivatives of the linear operators~\cite{piggyback2003,bogensperger22,chambolle2021learning}. This bears the advantage of not being limited to the number of primal-dual iterations due to computational memory issues. While an estimate for a saddle point for \eqref{eq:sp_general} is obtained, the adjoint state of the corresponding bi-quadratic saddle point problem is simultaneously computed (see Algorithm~\ref{alg:piggyback}). Using the resulting approximate saddle point $(u^J,v_K^J,v_L^J,p^J)$ and its adjoint state $(U^J,V_K^J,V_L^J,P^J)$, the gradients with respect to $\Kop$ and $\Lop$ can then be computed using automatic differentiation (see~\cite{chambolle2021learning} for more details).
\begin{algorithm}[ht]
\caption{Piggyback primal-dual algorithm to solve \eqref{eq:sp_general} and its adjoint.}
\SetAlgoLined 
\KwIn{initial values $(u^1,v_K^1,v_L^1,p^1)$ and $(U^1,V_K^1,V_L^1,P^1)$, block-preconditioned step size parameters  $\tau_u, \tau_{v_K},\tau_{v_L},\sigma>0$, $\theta\in[0,1]$, number of iterations $J$}
\KwResult{approximate saddle point $(u^J, v_K^J, v_L^J, p^J)$ and its adjoint state $(U^J, V_K^J, V_L^J, P^J)$}
\For{$j = 1, 2,\ \ldots ,\  J $}
{
$p^{j+1} =p^j + \sigma (\Dtwo u^j - \Epsop \Lop^\ast v_L^j - \Kop^\ast v_K^j), \quad P^{j+1}=P^j + \sigma (\Dtwo U^j - \Epsop \Lop^\ast V_L^j - \Kop^\ast V_K^j)$\;

$\bar{p}^{j+1} = p^{j+1} + \theta( p^{j+1} - p^j), \quad\qquad\qquad~\mspace{-3mu}\bar{P}^{j+1} = P^{j+1} + \theta( P^{j+1} - P^j)$\;

$\widetilde{u}^{j+1} = u^j - \tau_u D^{2^\ast} \bar{p}^{j+1}, \qquad\qquad\qquad\quad~\mspace{-1mu}\widetilde{U}^{j+1} = U^j - \tau_u (D^{2^\ast} \bar{P}^{j+1} + \nabla \ell(u^j,t)) $\;
$u^{j+1} = \prox_{\tau_u G(\cdot,f)}(\widetilde{u}^{j+1}), \qquad\quad\qquad\quad~\mspace{0mu} U^{j+1} = \nabla \prox_{\tau_u G(\cdot,f)}(\widetilde{u}^{j+1}) \cdot \widetilde{U}^{j+1}$\;

$\widetilde{v}_L^{j+1} =v_L^j + \tau_{v_L} \Lop \Epsop^\ast \bar{p}^{j+1}, \qquad\qquad\qquad\quad \mspace{-3mu}\widetilde{V}_L^{j+1}= V_L^j + \tau_{v_L} \Lop \Epsop^\ast \bar{P}^{j+1}$\;
 $v_L^{j+1}  = \prox_{\tau_{v_L} \Vert \cdot \Vert_{Z}} (\widetilde{v}_L^{j+1}),\qquad\qquad\qquad\quad \mspace{-3mu}V_L^{j+1}  = \nabla \prox_{\tau_{v_L \Vert \cdot \Vert_{Z}}} (\widetilde{v}_L^{j+1}) \cdot \widetilde{V}_L^{j+1}$ \;

$\widetilde{v}_K^{j+1} = v_K^j + \tau_{v_K} \Kop \bar{p}^{j+1},\qquad\qquad\qquad\qquad\mspace{-6mu} \widetilde{V}_K^{j+1} = V_K^j + \tau_{v_K} \Kop \bar{P}^{j+1} 
$\;
$v_K^{j+1} = \prox_{\tau_{v_K} \Vert \cdot \Vert_{Z}} (\widetilde{v}_K^{j+1}), \qquad\qquad\qquad\quad \mspace{-4mu}V_K^{j+1} = \nabla \prox_{\tau_{v_K} \Vert \cdot \Vert_{Z}} (\widetilde{v}_K^{j+1}) \cdot \widetilde{V}_K^{j+1} $\; }
\label{alg:piggyback} 
\end{algorithm} 

The outer bilevel learning problem is solved using a block-wise Adam optimizer~\cite{kingma2014adam}, whose block-wise structure is crucial due to the imposed constraints on the filters for the projections. 
This allows for individual adaptive learning rates for all groups of parameters subject to the same constraint by estimating the first and second gradient moments.

\subsection{Filter Settings}\label{sec:filtersettings}
\subsubsection{Initialization} Since the underlying problem is of a non-convex nature due to its bilevel structure we lack any guarantee to obtain a global minimum. Initialization can thus make a huge difference. Experiments with different initialization schemes were conducted, comparing filter coefficients drawn from a uniform or normal distribution, using the recently proposed discretization~\cite{hosseini2022second}, or using reference-style filters that introduce no sort of initial interpolation. We empirically found the initialization from~\cite{hosseini2022second} to work well for small filter kernel sizes of $3\times 3$ and $n_K,n_L \leq 4$, whereas for larger filter kernels uniformly distributed filters $\sim \mathcal{U}(-1/\sqrt{b},1/\sqrt{b})$ ($b$ depends on the number of input dimensions and the filter kernel size) yielded the most satisfactory results, which is inspired by the well-known Xavier initialization~\cite{glorot2010understanding}.

\subsubsection{Constraints} 
As given in Section~\ref{sec:learning}, the constraints are used to ensure the boundedness of the filter coefficients, such that the sum of each filter is constrained to be 1, i.e. 
\begin{equation}
    \sum_{m,n} \xi^{1,r}_{m,n} = \sum_{m,n}  \xi^{2,r}_{m,n} = \sum_{m,n}  \xi^{3,r}_{m,n} = 1, \notag\qquad
    \sum_{m,n} \eta^{1,l}_{m,n} = \sum_{m,n}  \eta^{2,l}_{m,n} = 1,
\end{equation}
with $r=1,\dots, n_K$, $l=1,\dots,n_L$. The corresponding projection per filter is computed following~\cite{chambolle2021learning}. The fact that second-order \gls{tgv} requires choosing two hyperparameters $\alpha_0$ and $\alpha_1$ majorly influencing the resulting reconstructions, where proper tuning can be challenging especially for natural images. Therefore an option, in this case, is to implicitly include them in the aforementioned constraints of the learned filters, such that the filter coefficients sum up to the same values $\gamma_K,\gamma_L \in \R$, respectively.

Moreover, one can also attempt to include a symmetry constraint to construct filters with a 90$^\circ$ rotational invariance property on the filter coefficients of $\Lop$~\cite{chambolle2021learning}. This reduces the actual number of learnable filters, which can be seen as an additional form of regularization in the learning setting to reduce overfitting to the training data. 
\section{Numerical Results}
\subsection{Data Sets} \label{sec:datasets} 
\subsubsection{Synthetic Data} 
\begin{figure}
    \centering
    \includegraphics[width=\textwidth]{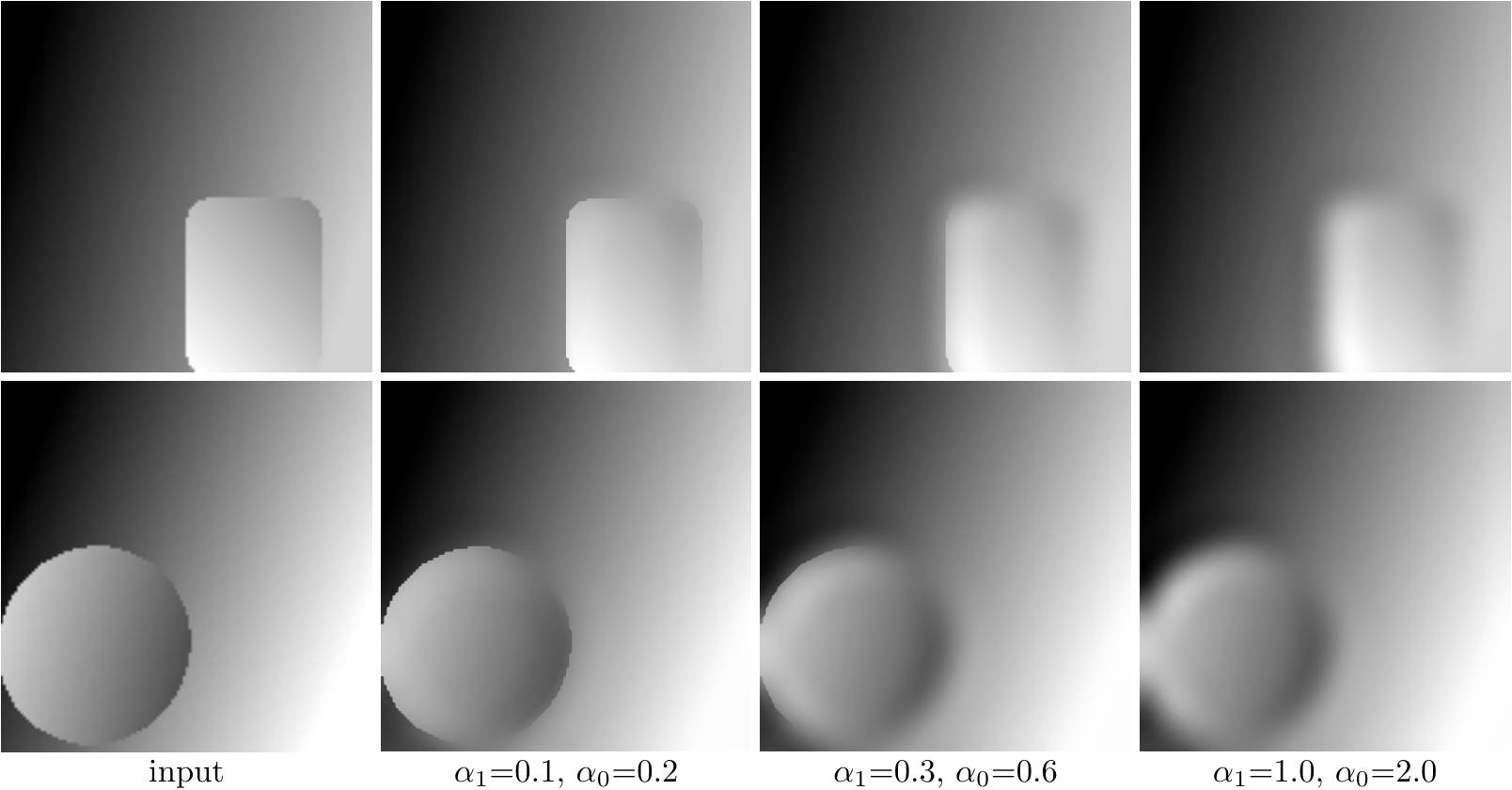}
    \caption{Sample input images (left column) and respective ``ground truth'' reconstructions obtained from applying Condat's inspired \gls{tgv} on upsampled images from the synthetic data set. This is done for three different combinations of $(\alpha_1,\alpha_0)$.}
    \label{fig:synthetic_images}
\end{figure}

A synthetic data set was generated which is inspired by the intrinsic nature of the second-order \gls{tgv} regularizer that favors piecewise affine solutions. A train and test data set each with 32 images of size $128 \times 128$ was constructed by randomly drawing basic shapes such as triangles, rectangles, and circles with varying sizes, which were filled with piecewise affine intensity changes and embedded within different (piecewise affine) background scenes. Examples of such images can be seen in the first column in Figure~\ref{fig:synthetic_images}. Casting this as an inverse problem requires some sort of ground truth to compare the obtained reconstructions for specific $(\alpha_1,\alpha_0)$. Although there exist special cases such as specific 1D functions where an actual solution for 
$\text{TGV}_\alpha^2$ exists~\cite{poschl2013exact}, there is no ground truth available for arbitrary 2D images. Thus, the idea is to upsample the images (we use a size of $8M \times 8N$) and to compute a pseudo ground truth using the Condat-inspired \gls{tgv}~\cite{hosseini2022second} due to its rotational invariance, where the intuition is that this solution better approximates the ground truth. 
A downsampled version of this is used as a new ground truth to compare the effect of using different handcrafted and learned discretization versions of $\text{TGV}_\alpha^2$. To enable a fair comparison, this was done for three distinct parameter settings $(\alpha_1, \alpha_0) \in \{(0.1,0.2), (0.3,0.6), (1.0,2.0)\}$ leading to different levels of smoothing, which is shown in the last three columns in Figure~\ref{fig:synthetic_images}. 
\subsubsection{Natural Images} Furthermore, a second, distinct data set was used. It is comprised of natural images where images were sampled from the well-known BSDS500 data set~\cite{arbelaez2010contour}. A train and a test data set each containing 32 images of size $128 \times 128$ were generated and all images were corrupted using zero-mean additive Gaussian noise $\sim \mathcal{N}(0,\sigma^2)$, which was sampled independently per pixel using noise levels $\sigma \in \{12.75,25.5\}$ corresponding to 5$\%$ or $10\%$ noise, respectively.

\subsection{Results}
\subsubsection{Synthetic Data}
The results presented here compare the solutions for the synthetic test data set using a ``denoising'' data term for different discretization schemes to the computed pseudo ground truth. In essence, the standard \gls{tgv} and two scenarios of each handcrafted and learned discretizations are compared, where for both the two settings $n_K=1$, $n_L=3$ as in~\cite{hosseini2022second} and $n_K=4$, $n_L=4$ are used. For the learning setting mostly $3000-5000$ iterations in the outer learning problem were used or until the training loss function had stabilized, while the inner problem was solved with 100 piggyback primal-dual iterations using a warm-starting initialization scheme in each learning step. For the evaluation with Algorithm~\ref{alg:pd_general} the number of iterations was substantially increased to ensure absolute convergence. 

Quantitative results reporting the mean \gls{psnr} and \gls{mse} on the test set are presented in Table~\ref{tab:exacttgv}. They clearly show
that both the handcrafted and the learned discretizations outperform the standard \gls{tgv}, which is reflected for all three scenarios of $(\alpha_1,\alpha_0)$. In each case using handcrafted filters gives some improvement, however directly learning the filters always outperforms these results by a larger margin. For $(\alpha_1,\alpha_0)=(1.0,2.0)$, which introduces a substantial amount of smoothing and is, therefore, more challenging, the increase due to the handcrafted filters is barely present, while the learned filters still manage to improve the result further. 
\begin{table}
\begin{center}
\caption{Quantitative comparison for the synthetic test images of the standard \gls{tgv} and different handcrafted and learned discretizations.}\label{tab:exacttgv}
\resizebox{\textwidth}{!}{
\begin{tabular}{p{4.5cm}|p{1.4cm} p{1.5cm}|p{1.4cm} p{1.5cm}|p{1.4cm} p{1.5cm}}
& \multicolumn{2}{c}{$\alpha_1$=0.1, $\alpha_0$=0.2} & \multicolumn{2}{c}{$\alpha_1$=0.3, $\alpha_0$=0.6}  & \multicolumn{2}{c}{$\alpha_1$=1.0, $\alpha_0$=2.0}  \\ \cline{2-7}
\backslashbox[49mm]
{\footnotesize{\textbf{Method}}}{\footnotesize{\textbf{Metric}}}
 &  PSNR & MSE $\cdot 10^{-2}$ & PSNR & MSE $\cdot 10^{-2}$ & PSNR & MSE $\cdot 10^{-2}$ \\
\hline\hline
\gls{tgv} &40.36 & 0.0135 &39.60 &0.0151  &35.36&0.0302 \\
Handcrafted Disc. $n_K$=1, $n_L$=3&  41.11& 0.011 &40.09  &0.0128&35.43  &0.0297  \\
Handcrafted Disc. $n_K$=4, $n_L$=4& 41.11& 0.011& 40.09 & 0.0128 &35.42 &0.0297  \\
Learned Disc. $n_K$=1, $n_L$=3, $3\times 3$ &43.05  &0.0065  & 41.37 &0.0088 &36.84&0.0218\\
Learned Disc. $n_K$=4, $n_L$=4, $3\times 3$ &\textbf{43.09} &\textbf{0.0064}  & \textbf{41.45} &\textbf{0.0086} &\textbf{36.95} &\textbf{0.0211}\\
\hline
\end{tabular}}
\end{center}
\end{table} 

\subsubsection{Natural Images}
Moreover, results on natural image denoising for 5\% and 10\% additive Gaussian noise are presented, again comparing different handcrafted and learned discretizations. As for the previous task, the number of filters $n_K$ and $n_L$ is varied, however, due to the nature of natural images it is reasonable to also use a higher number of filters and larger kernel sizes. Thus, we extend the kernel size to $7 \times 7$, which is a good trade-off in terms of increased globality while maintaining moderate complexity. It is noteworthy that using the largest filter settings results in a computational time increase of up to twelve-fold, however, this trade-off is justified by the clear performance improvements.
The learning setting remains similar to the previous experiments, and evaluations were conducted using $10^{4}$ primal-dual iterations. The hyperparameters were set to $\alpha_1=\{0.03,0.0685\}$ per noise level, respectively, resulting from a prior grid search (and $\alpha_0=2\alpha_1$). Note that due to the constraint that the filter coefficients are allowed to sum up to $\gamma_K$, $\gamma_L$ for the learning setting, the values of $\alpha_0$, $\alpha_1$ serve as an initialization, while their final values amount to $|\gamma_K| \alpha_0$ and $|\gamma_L| \alpha_1$ (whilst normalizing the learned filters with $\gamma_K$, $\gamma_L$). 

Quantitative results are summarized in Table~\ref{tab:natural_v2}, clearly showing that learned discretizations with a higher number of filters (such as $n_K$=16 and $n_L$=16) and a filter kernel size of $7\times7$ yield the best results in terms of \gls{psnr} in dB and \gls{mse} (improvements up to approx. 0.6~dB). 
The same is confirmed when evaluating the \gls{ssim}~\cite{ssim2004}. 
This can be expected as these settings allow us to learn a more complex and rich discretization of natural images.
The additional symmetry constraint on $\Lop$ -- indicated by (sym.) -- does not influence the results significantly. Further, an additional quantitative comparison using a \gls{tv} regularizer with hand-tuned $\alpha = \{0.03,0.0685\}$ for both noise levels confirms the well-established fact that the \gls{tv} is a proper handcrafted regularizer especially considering its simplicity. 

\begin{table}
\begin{center}
\caption{Quantitative comparison of natural image denoising of the test set with 5\% and 10\% Gaussian noise for different handcrafted and learned discretizations.}\label{tab:natural_v2}
\resizebox{\textwidth}{!}{
\begin{tabular}{p{5.7cm}|p{1.1cm} p{1.1cm} p{1.1cm} |p{1.1cm} p{1.1cm} p{1.1cm}}
\multirow{2}{*}{ \backslashbox[60mm]{\footnotesize{\textbf{Method}}}{\footnotesize{\textbf{Metric}}}}  &  \multicolumn{3}{c|}{5\% Gaussian noise}  & \multicolumn{3}{c}{10\% Gaussian noise} \\ \cline{2-7}
 &  PSNR & MSE $\cdot 10^{-2}$ & SSIM & PSNR & MSE$\cdot 10^{-2}$ & SSIM \\
\hline\hline
Corrupted $f$& 26.04 & 0.2490 & 0.7885 &  20.02& 0.9959 & 0.5382 \\
\gls{tv} & 30.14& 0.1049&0.9249 & 26.52 & 0.2445 & 0.8497 \\
\gls{tgv} & 30.2 & 0.1043& 0.9257 &26.56 &0.2431 & 0.8512 \\
Handcrafted Disc. $n_K$=1, $n_L$=3& 30.24 & 0.1046 & 0.9267 & 26.69&0.2394 & 0.8553 \\
Handcrafted Disc. $n_K$=4, $n_L$=4& 30.29 & 0.1030 &0.9278 &26.71 &0.2370 & 0.8565 \\
Learned Disc. $n_K$=1, $n_L$=3, $3\times 3$ & 30.52 & 0.0935 & 0.9274 &  26.95 & 0.2172 & 0.8596\\
Learned Disc. $n_K$=4, $n_L$=4, $3\times 3$ &30.66 &0.0906&0.9298 & 27.06  &0.2123& 0.8620\\
Learned Disc. $n_K$=8, $n_L$=8, $7\times 7$ & 30.74 & 0.0896& 0.9314& 27.14 &0.2090& 0.8649\\
Learned Disc. $n_K$=8, $n_L$=8, $7\times 7$, sym. &30.72 &0.0898 &0.9311&27.15 &0.2089 & 0.8649 \\
Learned Disc. $n_K$=10, $n_L$=10, $7\times 7$ &30.73 & 0.0896 & 0.9313 & 27.17 & 0.2081 & 0.8657\\
Learned Disc. $n_K$=16, $n_L$=16, $7\times 7$  & 30.77 &0.0891& 0.9319 & 27.16& 0.2087 & 0.8654\\
Learned Disc. $n_K$=16, $n_L$=16, $7\times 7$, sym. & \textbf{30.77} & \textbf{0.0890} & \textbf{0.9320} & \textbf{27.18}& \textbf{0.2074}& \textbf{0.8659} \\
\hline
\end{tabular}}
\end{center}
\end{table} 

\begin{figure}[H]
    \centering
    \includegraphics[width=\textwidth]{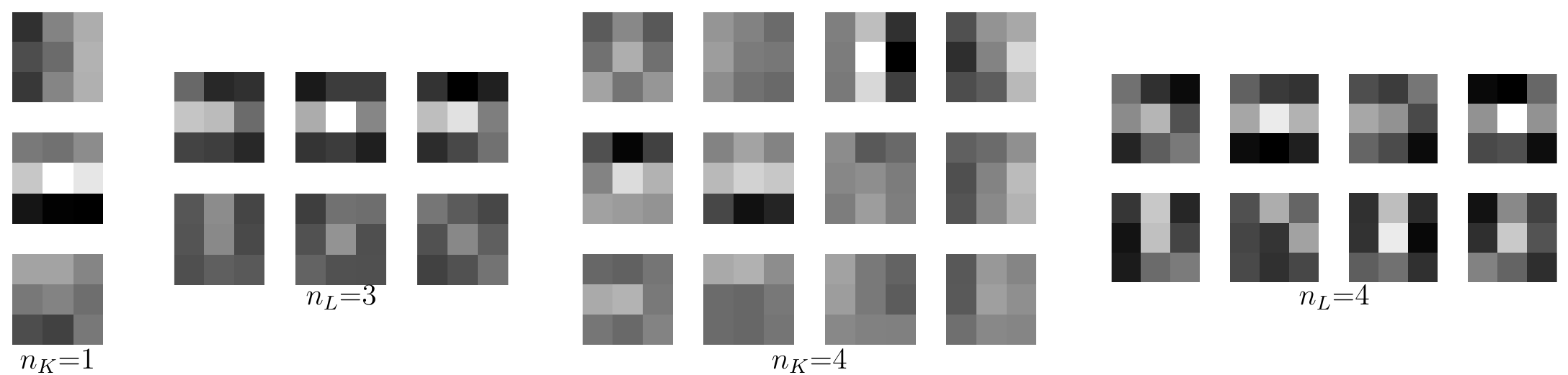}
    \caption{Learned interpolation filters for the setting $n_K$=1,$n_L$=3 and $n_K$=4,$n_L$=4, also cf. Figure~\ref{fig:handcrafted} for a visual comparison with the handcrafted filters.}
    \label{fig:learned_small}
\end{figure} 

\begin{figure}[H]
    \centering
    \includegraphics[width=\textwidth]{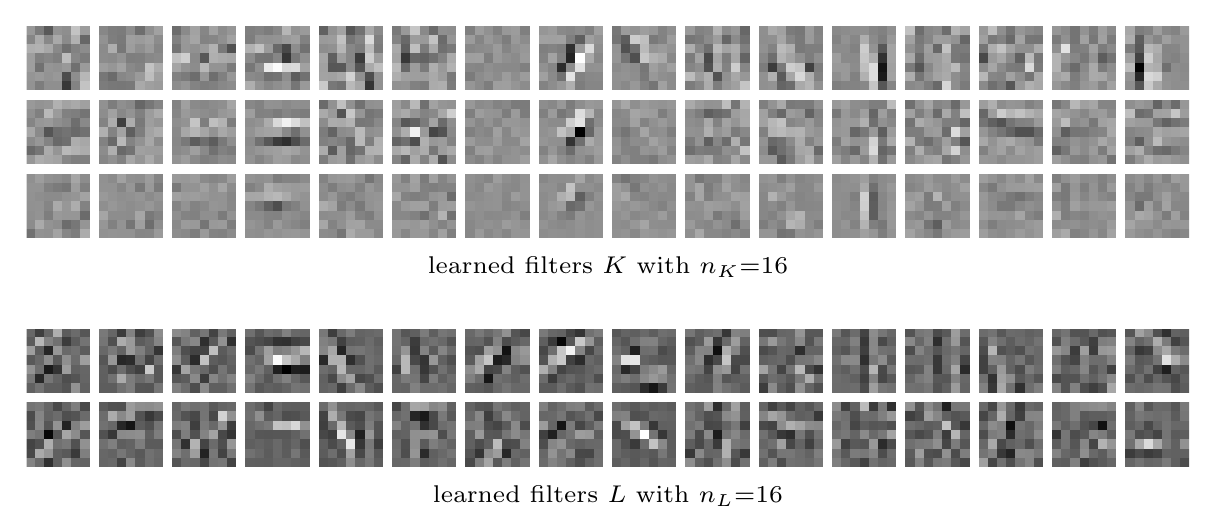}
    \caption{Learned $7\times 7$ filters  using $n_L$=16 and $n_K$=16 for denoising (10\% Gaussian noise). The row of a depicted filter denotes the component of the respective vector/tensor field that it acts upon, whereas the column refers to the specific filter $r$ or $l$ (with $r=1,\cdots,n_K$ and $l=1,\cdots,n_L$.)}
    \label{fig:naturalfilters} 
\end{figure} 

\begin{figure}[H]
    \centering
    \includegraphics[width=\textwidth]{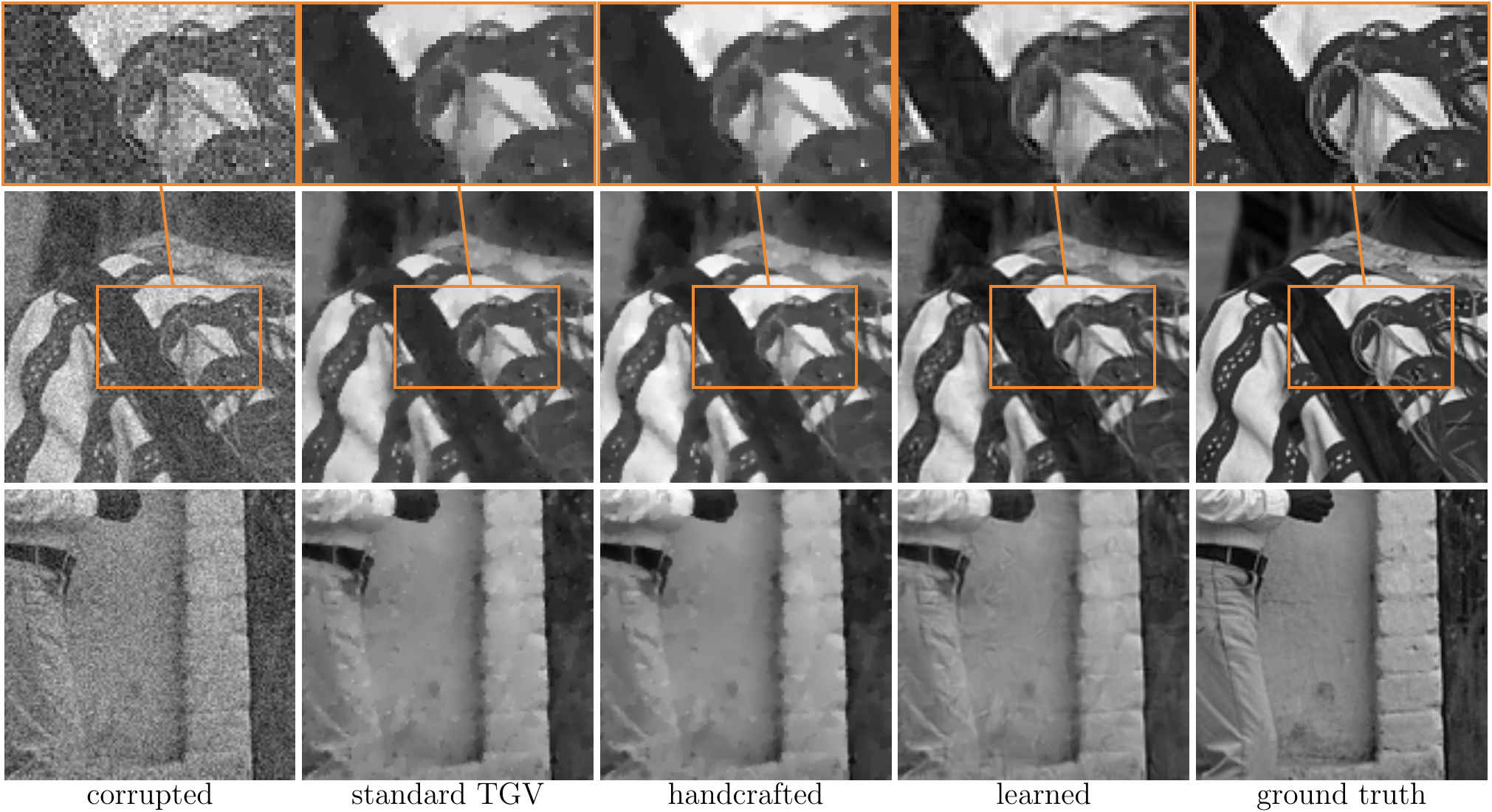}
    \caption{Sample reconstructions from natural test images (10\% Gaussian noise) comparing the standard \gls{tgv}, the handcrafted discretization scheme with $n_K$=1, $n_L$=3~\cite{hosseini2022second}, and learned filters using $n_K$=16, $n_L$=16. For completeness, the ground truth images are also shown.}
    \label{fig:natural} 
\end{figure}

Exemplary learned filters $\Kop$ and $\Lop$ are displayed in Figure~\ref{fig:learned_small} for filter settings as shown in Figure~\ref{fig:handcrafted} and further in Figure~\ref{fig:naturalfilters} for larger filters using the settings $n_K=16$ and $n_L=16$. 
Generally, it can be observed that different orientations are captured, while some of the filters in $\Kop$ introduce a bit of a smoothing effect. This can be associated with the fact that this filter operates on second-order finite difference arrays that are already very smooth, thus the discretization at this level will not contribute as much as opposed to the filters contained in $\Lop$. Qualitative results on 10\% Gaussian noise image denoising are shown in Figure~\ref{fig:natural}, where the handcrafted filters from~\cite{hosseini2022second} and the learned filters with $n_K$=16 and $n_L$=16 are compared in terms of reconstruction quality of two sample test images. Using learned filters tends to exhibit finer details and produces significantly more structure in the reconstructed images. 

\section{Conclusion}
We proposed a general discretization scheme for the second-order \gls{tgv} regularizer building upon the idea of~\cite{hosseini2022second}. This is supported by a proof of consistency by means of $\Gamma$-convergence of the newly discretized functional. Moreover, using a synthetic and a natural image data set we showcase that learning interpolation filters quantitatively and qualitatively improves the resulting image reconstruction in the setting of image denoising. It suggests that there might not be an ideal predefined set of interpolation filters applicable to all data sets and image reconstruction settings, but the most suited one can be obtained by learning within the respective setting. 
The proposed framework can be adapted to higher-order versions of \gls{tgv} and further to other linear inverse problems. This is subject to future work, as well as an analysis on the generalization of the learned filters and on the robustness in terms of rotational invariance. 

\paragraph{Acknowledgements} Lea Bogensperger and Thomas Pock acknowledge support by the BioTechMed Graz flagship  project ``MIDAS''. Alexander Effland was funded by the German Research Foundation under Germany’s Excellence Strategy – EXC-2047/1 – 390685813 and – EXC2151 – 390873048.

\appendix
\setcounter{equation}{0}
\renewcommand{\theequation}{A.\arabic{equation}}
\section{Appendix: Consistency}
For consistency, we have to show the $\Gamma$-convergence of the discretized second-order \gls{tgv} to the continuous definition, where the latter reads as
\begin{align}\label{eq:tgv_cont}
    \text{TGV}_{\alpha}^2(u) & = \sup_{p} \bigg \{ \int_\Omega u ~\Divtwo p \d{x}: p \in \mathcal{C}^\infty(\Omega, \text{Sym}^{2\times 2}), \Vert p \Vert_{\infty} \leq \alpha_0, \Vert \Div p \Vert_{\infty} \leq \alpha_1 \bigg \} \notag \\ &= 
    \min_{w \in BD(\Omega)} \bigg \{ \int_\Omega \alpha_1 |\nabla u - w | + \alpha_0 |\Epsop w | \d{x} \bigg \}.
\end{align}
The latter variational problem turns out to always have
a minimizer $w$, thanks to a Rellich type compactness theorem in $BD(\Omega)$ (the space of fields such that
$Ew$ is a bounded measure), see~\cite {TemamStrang}.

Consider the domain $\Omega=(0,1)^2$ of $N\times N$ pixels, where each pixel is of size $h\times h$. The operators and variables are marked with an $h$ to denote that the grid is discretized in steps of $h$. For $u^h \in \domain$, $w^h \in \domain^2$, $v_K^h \in \domain^{3\times n_K}$, $v_L^h \in \domain^{2\times n_L}$, $p^h \in \domain^3$ we use both primal and dual definitions of the discretized second-order \gls{tgv}
\begin{align} \label{eq:tgv_disc}
    \text{TGV}^2_{\alpha,h}(u^h) & =
    \sup_{p^h}\bigg\{h^2 \langle \Div_h^2 p^h, u^h \rangle : \Vert \Lop_h \Div_h p^h \Vert_Z^\ast \leq \alpha_1, \Vert \Kop_h p^h \Vert_Z^\ast \leq \alpha_0 \bigg \}
    \\ &=  
    \min_{w^h,v_K^h,v_L^h} \bigg\{ h^2 \alpha_1 \Vert v_L^h\Vert_{Z} + h^2\alpha_0 \Vert v_K^h \Vert_{Z}: \Lop_h^\ast v_L^h = \D_h u^h - w^h \text{, } \Kop_h^\ast v_K^h = \Epsop_h w^h \bigg \}. \notag
\end{align}
In a slightly simpler setting where we assume that $u$ is global affine plus periodic and $w$ is periodic with periodic boundary conditions, theorem~\ref{th:consistency} holds, where the corresponding proof closely follows the respective proof in~\cite{chambolle2021learning}.
\begin{proof}
For the $\Gamma$-lower limit, we consider an image $u$ and a sequence of discrete images $u^h$ which, viewed as piecewise constant functions on pixels of size $h\times h$, converge to $u$ in $L^1(\Omega)$ as $h\to 0$.
Then, we consider a dual test field $p \in \mathcal{C}^\infty(\Omega, \text{Sym}^{2\times 2})$ with compact support satisfying the constraints in~\eqref{eq:tgv_cont}. 
We have to find a discretization of $p$ such that
\begin{equation*}
\int_\Omega u  \Divtwo p \d{x} = \lim_{h\to 0} h^2 \langle u^h, \Div_h^2 p^h \rangle.
\end{equation*}
Since $p$ is smooth we can simply consider its discretization $p_{i,j}^{h,\bullet}=\bbar{p}^\bullet(ih,jh)$ inside the pixel $c_{ij}=\big((i-\tfrac 1 2)h,(i+\tfrac 1 2)h)\big)\times \big((j-\tfrac 1 2)h,(j+\tfrac 1 2)h)\big)$. 
As $u^h$ converges to $u$ in $L^1$, we have to prove the uniform convergence of $p^h$ and its discrete derivatives to $p$ and its corresponding derivatives, where we in particular have to show
\begin{equation}\label{eq:divcont_eq_divdisc}
\Div^2 p (ih,jh) \approx \Div_h^2 p^h
\end{equation}
up to some errors which uniformly converge to $0$ as $h\to 0$.
In the continuous setting of the left side, this can be expressed as
\begin{equation}
    \div^2 p = \div \begin{pmatrix} \partial_1 p^1  + \partial_2 p^2 \\ \partial_1 p^2 + \partial_2 p^3 \end{pmatrix} = \partial_1 \partial_1 p^2 + \partial_1 \partial_2 p^2 + \partial_2 \partial_1 p^2 + \partial_2 \partial_2 p^3. \notag
\end{equation}
In the discrete setting, the right side of~\eqref{eq:divcont_eq_divdisc} can be written as
\begin{align}\label{eq:a}
    (\div_h (\div_h p^h))_{i,j} &= \tfrac 1 h \div_h \begin{pmatrix}
        p^{h,1}_{i+1,j}-p^{h,1}_{i,j} + p^{h,2}_{i,j+1}-p^{h,2}_{i,j} \\ p^{h,2}_{i+1,j} - p^{h,2}_{i,j} + p^{h,3}_{i,j+1} - p^{h,3}_{i,j} \end{pmatrix} \\ 
        & = \tfrac{1}{h^2} \Big ( p^{h,1}_{i+2,j} - 2p^{h,1}_{i+1,j} +p^{h,1}_{i,j} + 2 p^{h,2}_{i+1,j+1} - 2p^{h,2}_{i+1,j} - 2 p^{h,2}_{i,j+1} \notag \\ 
        & \quad + 2p^{h,2}_{i,j} + p^{h,3}_{i,j+2} -2p^{h,3}_{i,j+1} + p^{h,3}_{i,j} \Big ). \notag
\end{align}
Since each individual component of $p^{h,1},p^{h,2},p^{h,3}$ is equal to the sampling of the corresponding component of $p$ at $(ih,jh)$ we can write
\begin{align*}
    p^{h,1}_{i+2,j} - 2p^{h,1}_{i+1,j} +p^{h,1}_{i,j} &= \bbar{p}^1((i + 2)h,jh) - 2\bbar{p}^1((i+1)h,jh) + \bbar{p}^1(ih,jh), \\
    2 p^{h,2}_{i+1,j+1} - 2p^{h,2}_{i+1,j} - 2 p^{h,2}_{i,j+1} + 2p^{h,2}_{i,j} & = 2\bbar{p}^2((i+1)h, (j+1)h) -2\bbar{p}^2((i+1)h,jh) \notag \\  & \quad -2\bbar{p}^2(ih,(j+1)h) + 2\bbar{p}^2(ih,jh), \notag \\
    p^{h,3}_{i,j+2} -2p^{h,3}_{i,j+1} + p^{h,3}_{i,j} &= \bbar{p}^3(ih,(j+2)h) -2\bbar{p}^3(ih,(j+1)h) + \bbar{p}^3(ih,jh). \notag
\end{align*}
Now, since $p$ is smooth, a second-order Taylor expansion can be applied for each individual term, i.e., for the first component of the tensor field $\bbar{p}^1$ this reads as follows (the third term needs no expansion):
\begin{align*}
    \bbar{p}^1((i+2)h,jh) &= \bbar{p}^1(ih,jh) + 2h \partial_1 \bbar{p}^1(ih,jh) + 2 h^2 \partial_1 \partial_1 \bbar{p}^1(ih,jh) + \mathcal{O}(h^3), \\
    -2\bbar{p}^1((i+1)h,jh) &= -2\bbar{p}^1(ih,jh) -2 h \partial_1 \bbar{p}^1(ih,jh) - h^2 \partial_1 \partial_1 \bbar{p}^1(ih,jh) + \mathcal{O}(h^3), \notag \\
    \bbar{p}^1(ih,jh) &= \bbar{p}^1(ih,jh), \notag
\end{align*}
whose sum equals $h^2 \partial_1 \partial_1 \bbar{p}^1(ih,jh) + \mathcal{O}(h^3)$.
Similarly, for the components involving $\bbar{p}^3$ the same procedure can be performed, which yields $h^2 \partial_2 \partial_2 \bbar{p}^3(ih,jh) + \mathcal{O}(h^3)$.
Here and in all that follows, $\mathcal{O}(h^3)$ only depends on a global bound on the third derivatives of the smooth field $p$.

A similar computation is performed for the mixed derivatives of the component $\bbar{p}^2$. 
Using a Taylor expansion for three of the terms we obtain
\begin{align*}
    -2 \bbar{p}^2((i+1)h,jh) &= -2\bbar{p}^2(ih,jh) - 2h \partial_1 \bbar{p}^2(ih,jh) - h^2 \partial_1 \partial_1 \bbar{p}^2(ih,jh) + \mathcal{O}(h^3), \\
    - 2\bbar{p}^2(ih,(j+1)h) &= -2\bbar{p}^2(ih,jh) -2 h \partial_2 \bbar{p}^2(ih,jh) - h^2 \partial_2 \partial_2 \bbar{p}^2(ih,jh) +\mathcal{O}(h^3), \notag \\
    2\bbar{p}^2((i+1)h,(j+1)h) &=2 \bbar{p}^2(ih,jh) +2 h \partial_1 \bbar{p}^2(ih,jh) + 2h \partial_2 \bbar{p}^2(ih,jh) + h^2 \partial_1 \partial_1 \bbar{p}^2(ih,jh) \notag \\ & \quad + h^2 \partial_2 \partial_2 \bbar{p}^2(ih,jh)  +  h^2 \partial_1 \partial_2 \bbar{p}^2(ih,jh) + h^2 \partial_2 \partial_1 \bbar{p}^2(ih,jh) + \mathcal{O}(h^3), \notag\\
    2\bbar{p}^2(ih,jh) &= 2\bbar{p}^2(ih,jh), \notag
\end{align*}
where only $h^2 \partial_1 \partial_2 \bbar{p}^2(ih,jh) + h^2 \partial_2 \partial_1 \bbar{p}^2(ih,jh) + \mathcal{O}(h^3)$ remains. 
Finally, this leads to (compare~\eqref{eq:a})
\begin{align*}
    \div_h^2 p_{i,j}^h =  \partial_1 \partial_1 \bbar{p}^1(ih,jh) + \partial_1 \partial_2 \bbar{p}^2(ih,jh)  + \partial_2 \partial_1 \bbar{p}^2(ih,jh) + \partial_2 \partial_2 \bbar{p}^3 (ih,jh)  + \mathcal{O}(h),
\end{align*}
which readily implies
\begin{equation*}
        \int_\Omega u \div^2 p \d{x} = h^2 \langle \div_h^2 p^h, u^h \rangle + \mathcal{O}(h).
\end{equation*}
It remains to slightly modify $p^h$ to satisfy the constraints of the discrete second-order \gls{tgv} in~\eqref{eq:tgv_disc}. 
We have to show that our discretization is admissible up to a small error $C h$, since the filters are constrained to sum to 1. 
Let us recall that $\nu \in \mathbb{N}$ denotes the support $(2\nu +1)\times(2\nu + 1)$ of the filter kernels.
For simplicity, we only consider the case $n_K=1$.
Hence, 
\begin{equation*}
    K_h^{1}p^{h,1}_{i,j} = \sum_{m,n=-\nu}^\nu \xi^{1}_{m,n} p^{h,1}_{i-m,j-n} = p^{1}(ih,jh) + \sum_{m,n=-\nu}^\nu \xi^{1}_{m,n} (p^{h,1}_{i-m,j-n} - p^1 (ih,jh)).
\end{equation*}
Since $p$ is smooth, there is a constant $C_1$ for $\Kop_h$ depending on $\Vert \nabla p\Vert_\infty$, such that the last term can be bounded by $\alpha_0 C_1 h/2$ (recall that we have a symmetric $2\times 2$ tensor field).
This implies that the error is of order $h$, therefore $\Vert \Kop_h p^h\Vert_Z^* \leq \alpha_0 ( 1+ C_1 h)$ and thus $\frac{p^h}{1+C_1 h}$ yields an admissible dual variable. 

One can proceed analogously for the filter $\Lop_h$, where we additionally have to incorporate the divergence operator in the constraint $\Vert \Lop_h \Div_h p^h \Vert_Z^\ast$.
Again, let us assume $n_L=1$.
Since $L_h^{1} (\Div_h p^h)^1=L_h^{1} (\partial_1^h p^{h,1} + \partial_2^h p^{h,2})$, we can examine the individual terms (the filter coefficients of $L_h^{l}$ are bounded and only have small support)
\begin{align*}
    L_h^{1} \partial_1^h p^{h,1} &= \sum_{m,n=-\nu}^\nu \eta_{m,n}^{1} \frac{p^{h,1}_{i+1-m,j-n} -p^{h,1}_{i-m,j-n} }{h} \\ &=\tfrac 1 h \sum_{m,n=-\nu}^\nu \eta_{m,n}^{1}\left(\bbar{p}^{1}(ih-mh+h,jh - nh) - \bbar{p}^{1}(ih-mh,jh - nh) \right).
\end{align*}
 Using the fundamental theorem of calculus, this yields
\begin{align*}
    &\tfrac 1 h \sum_{m,n=-\nu}^\nu \eta_{m,n}^{1}  \int_0^h \partial_1 \bbar{p}^{1}(ih-mh+t,jh-nh) \d{t}  \\ =&
    \partial_1 p^1(ih,jh) +  \tfrac 1 h \sum_{m,n=-\nu}^\nu \eta_{m,n}^{1} \int_0^h \big (\partial_1 \bbar{p}^{1}(ih-mh+t,jh-nh) - \partial_1 p^1(ih,jh) \big ). 
\end{align*}
As before, there exists a constant $C_2$ only depending on $\|\nabla^2 p\|_\infty$ and $h$ such that the last term can be bounded by $\frac{\alpha_1 C_2 h}{2\sqrt{2}}$, which implies
\begin{equation*}
    \Vert \Div_h p^h \Vert_Z^* \leq \Vert \Div p (ih,jh) \Vert_Z^* + \alpha_1 C_2 h.
\end{equation*}
Therefore, $\Vert \Div_h p^h \Vert_Z^* \leq \alpha_1 (1 + C_2 h)$ and thus setting $C=\max(C_1,C_2)$ we observe that $\frac{p^h}{1+C h}$ is an admissible dual variable.
Hence,
\begin{equation*}
h^2 \langle u^h, \Div_h^2 p^h \rangle\leq (1+C h ) \text{TGV}^2_{\alpha,h}(u^h).
\end{equation*}
Finally, letting $h\to 0$ and taking the supremum with respect to $p$ we obtain 
\begin{equation*}
    \mathrm{TGV}_{\alpha}^2(u) \leq \lim \inf_{h\to 0} \mathrm{TGV}^2_{\alpha,h} (u^h).
\end{equation*}

For the $\Gamma$-upper limit, we need to show that for any $u$ there exist discrete images $u^h$ such that $u^h\to u$ as $h\to 0$ and $\lim \sup_{h\to 0} \text{TGV}^2_{\alpha,h}(u^h)\le \text{TGV}_\alpha^2(u)$. 
We start from the primal definition in~\eqref{eq:tgv_cont} and only consider a $u$ which is of the form of a affine plus periodic function, and periodic $w$.
In this case, given $u$, if $u-q\cdot x$ is periodic, and $u^h$ is a recovery sequence for $u-q\cdot x$, then $(u^h_{i,j}+p\cdot (ih,jh))_{i,j}$ is a recovery sequence for $u$.
Without loss of generality, we can even assume that $u$ is periodic.
Given such a $u$, there exists at least one optimal $w$.

Let $\eta$ be a symmetric, non-negative mollifier and $\eta_\epsilon(x)=(1/\epsilon^2)\eta(x/\epsilon)$. 
We set
$u_\epsilon=\eta_\epsilon*u$ and $w_\epsilon=\eta_\epsilon*w$, and
by convexity we observe that
\[
\int_\Omega \alpha_1|Du_\epsilon-w_\epsilon|+\alpha_0|Ew_\epsilon|\le
\int_\Omega \alpha_1|Du-w|+\alpha_0|Ew|,
\]
so that
\[
\lim_{\epsilon\to 0}\int_\Omega \alpha_1|Du_\epsilon-w_\epsilon|+\alpha_0|Ew_\epsilon| =
\int_\Omega \alpha_1|Du-w|+\alpha_0|Ew| =\text{TGV}_\alpha^2(u),
\]
using the lower semicontinuity of these integrals.
Then, since $u_\epsilon$ and $w_\epsilon$ are smooth, we can approximate them by functions $u_{\epsilon,n}$ and $w_{\epsilon,n}$with a finite number of Fourier modes (by dropping all  Fourier coefficients with norm larger than $n\in\N$).
Hence,
\[
\lim_{n\to\infty}\int_\Omega \alpha_1|Du_{\epsilon,n}-w_{\epsilon,n}|+\alpha_0|Ew_{\epsilon,n}| =
\int_\Omega \alpha_1|Du_\epsilon-w_\epsilon|+\alpha_0|Ew_\epsilon|.
\]
As a result, by a standard diagonal argument, we can construct a sequence
$(u_k,w_k) = (u_{\epsilon_k,n_k},w_{\epsilon_k,n_k})$ which converges
to $(u,w)$ in $L^1$ and satisfies for $k \to \infty$
\begin{equation}\label{eq:approxfinitemodes}
    \int_\Omega \alpha_1 | D u_k - w_k | + \alpha_0 |\Epsop w_k | \d{x}  \to \text{TGV}_\alpha^2(u).
\end{equation}
For convenience of notation, we drop for a while the subscript $k$ and assume we are given $u,w$ with vanishing spectrum, i.e., $\hat{u}(\ell)=0$, $\hat{w}(\ell)=0$ if $|\ell | > R$. 
We can thus discretize for $N>R$ and $h=1/N$ as follows:
\begin{align*}
u_{i,j}^h = u(ih,jh) = \sum_{n,m=-N}^N \hat{u}_k(n,m) e^{i\pi \frac{in + jm}{N}}, \\
\left(w_{i+\frac12,j}^{h,1}, w_{i,j+\frac12}^{h,2} \right) = w(ih,jh) = \sum_{n,m=-N}^N \hat{w}_k(n,m) e^{i\pi \frac{in + jm}{N}}. \notag 
\end{align*}
The objective is to find $v_K^h$ and $v_L^h$ satisfying the constraints in the primal definition in~\eqref{eq:tgv_disc} such that
\begin{equation}\label{eq:limsup_general}
    \text{TGV}_{\alpha,h}^2 \big (u(ih,jh)\big )\leq h^2 \alpha_0 \Vert v_K^h \Vert_Z + h^2 \alpha_1 \Vert v_L^h \Vert_Z \lesssim \int_{\Omega} \alpha_1 | D u - w | + \alpha_0 |\Epsop w| \d{x} .
\end{equation}
The definitions of the discrete finite difference operators can be used to rephrase the objective to finding $v_K^h$ and $v_L^h$ such that
\begin{align}\label{eq:tgv_limsup_aim}
&h^2\alpha_1 \Vert v_L^h \Vert_Z + h^2 \alpha_0 \Vert v_K^h \Vert_Z \\ 
\lesssim &  h^2 \sum_{i,j=1}^N \alpha_1 \sqrt{\Big (\tfrac{1}{h}(u^h_{i+1,j} - u^h_{i,j})-  w^{h,1}_{i+\frac12,j}\Big)^2 + \Big (\tfrac{1}{h}(u^h_{i,j+1} - u^h_{i,j}) -  w^{h,2}_{i,j+\frac12}\Big )^2} + \notag  \\ 
&h^2  \sum_{i,j=1}^N \alpha_0  \sqrt{\tfrac{1}{h^2}( w^{h,1}_{i+\frac32,j}-w^{h,1}_{i+\frac12,j})^2 +\tfrac{1}{h^2}(w^{h,2}_{i,j+\frac32}-w^{h,2}_{i,j+\frac12})^2 \rule{0pt}{10pt}\ } \notag \\
& \quad \overline{\rule{0pt}{17pt}
    + \tfrac{2}{h^2} \Big (\tfrac 1 2 (w^{h,1}_{i+\frac12,j+1} - w^{h,1}_{i+\frac12,j}) + \tfrac 1 2 (w^{h,2}_{i+1,j+\frac12} - w^{h,2}_{i,j+\frac12})\Big )^2
\ }.
 \notag
\end{align}
First, let us find $v_L^h$ that fulfills the first term in~\eqref{eq:tgv_limsup_aim}. 
For simplicity, we consider only one filter component, i.e. $n_L=1$.
In essence, we want to find $v_L^{h,1}$ and $v_L^{h,2}$ such that
\begin{align}  \label{eq:Du_w=vL}
    \Big (\frac{u^h_{i+1,j}-u^h_{i,j}}{h} - w^{h,1}\Big )_{i+\frac12,j} = \sum_{m,n=-\nu}^\nu \eta^1_{m,n} v^{h,1}_{L_{i+m,j+n}}, \\
    \Big(\frac{u^h_{i,j+1}-u^h_{i,j}}{h}  - w^{h,2}\Big)_{i,j+\frac12} = \sum_{m,n=-\nu}^\nu \eta^2_{m,n} v^{h,2}_{L_{i+m,j+n}}. \notag
\end{align}
We now analyze in detail the first line in~\eqref{eq:Du_w=vL}, the computation for the second equation is analogous.
The discrete Fourier transform for $(r,s)\in \mathbb{Z}^2$ reads as
\begin{align}\label{eq:limsup_dft}
    &\frac{1}{2N}\sum_{i,j=-N}^N  \Big (\frac{u^h_{i+1,j}-u^h_{i,j}}{h} e^{-i\pi \frac{ir + js}{N}} - w^{h,1}_{i+\frac12,j} e^{-i\pi \frac{ir + js}{N}} \Big )\\
    =& \frac{1}{2N}\sum_{i,j=-N}^N \sum_{m,n=-\nu}^\nu \eta^1_{m,n} v^{h,1}_{L_{i+m,j+n}} e^{-i\pi \frac{ir + js}{N}}.\label{eq:limsup_dftsecondline}
\end{align}
Then, \eqref{eq:limsup_dft} can be written as
\begin{align*}
    & \frac{1}{2N}\sum_{i,j=-N}^N  \Big (\frac{u^h_{i+1,j}}{h}e^{-i\pi \frac{(i+1)r+js}{N}} e^{i\pi \frac{r}{N}} - \frac{u^h_{i,j}}{h}e^{-i\pi \frac{ir+js}{N}} - w^{h,1}_{i+\frac12,j}e^{-i\pi \frac{ir+js}{N}}\Big ) \\
    =& \hat{u}(r,s)\frac{e^{i\pi r h} - 1}{h} - \hat{w}(r,s). \notag
\end{align*}
Moreover, using that $\sum_{m,n=-\nu}^\nu \eta^1_{m,n} = 1$, \eqref{eq:limsup_dftsecondline} can be expressed as
\begin{align*}
    &\frac{1}{2N} \sum_{i,j=-N}^N \sum_{m,n=-\nu}^\nu \eta^{1}_{m,n} e^{i\pi \frac{mr +ns}{N}} v^{h,1}_{L_{i+m,j+n}} e^{-i\pi \frac{(i+m)r + (j+n)s}{N}} = \hat{v}^1_L(r,s) \sum_{m,n=-\nu}^\nu \eta^1_{m,n} e^{i\pi h (mr + ns)} \\
    =& \hat{v}^1_L(r,s) \Big ( 1 + \sum_{m,n=-\nu}^\nu \eta^1_{m,n} \big (e^{i\pi h (mr + ns)} - 1 \big ) \Big). 
\end{align*}
Thus, we obtain an expression for $\hat{v}^1_L$, which ensures $\hat{v}^1_L(r,s) = 0$ for $\Vert(r,s)\Vert > R$ since $\hat{u}$ and $\hat{w}$ only have a finite number of modes, which reads as
\begin{equation*}
    \hat{v}^1_L(r,s) = \frac{\hat{u}(r,s)\frac{e^{i\pi r h} - 1}{h} - \hat{w}(r,s)}{1 + \sum_{m,n=-\nu}^\nu \eta^1_{m,n} \big (e^{i\pi h (m r + n s)} - 1 \big )}.
\end{equation*}
To show the first part of~\eqref{eq:tgv_limsup_aim}, we use the inverse discrete Fourier transform to get
\begin{equation*}
    v^{h,1}_{L_{i,j}} = \sum_{r,s=-N}^N \hat{v}^1_L(r,s) e^{i\pi \frac{i r + j s}{N}} = \sum_{\Vert (r,s)\Vert \leq R}  \frac{\hat{u}(r,s)\frac{e^{i\pi r h} - 1}{h} - \hat{w}(r,s) }{1 + \sum_{m,n=-\nu}^\nu \eta^1_{m,n} \big (e^{i\pi h(m r + n s)} - 1 \big )} e^{i\pi \frac{i r + j s}{N}},
\end{equation*}
which can be used to express $\big(\frac{u^h_{i+1,j}-u^h_{i,j}}{h} - w^{h,1}\big )_{i+\frac12,j} - v^{h,1}_{L_{i,j}}$. Therefore we use
\begin{equation*}
\frac{1}{2N} \sum_{i,j=-N}^N \Big (\frac{u^h_{i+1,j}-u^h_{i,j}}{h} - w^{h,1}_{i+\frac12,j} \Big ) e^{-i\pi \frac{i r + j s }{N}} = \hat{v}^1_L(r,s) \sum_{m,n=-\nu}^\nu \eta^1_{m,n} e^{i\pi \frac{m r + n s}{N}},
\end{equation*}
such that we can obtain
\begin{equation*}
\frac{u^h_{i+1,j}-u^h_{i,j}}{h} - w^{h,1}_{i+\frac12,j} = \sum_{\Vert (r,s)\Vert \leq R} \hat{v}^1_L(r,s) \sum_{m,n=-\nu}^\nu \eta^1_{m,n} e^{i\pi \frac{m r + n s }{N}}e^{i\pi \frac{i r + j s}{N}}.
\end{equation*}
Finally, this yields
\begin{align*}
    &\Big(\frac{u^h_{i+1,j}-u^h_{i,j}}{h} - w^{h,1}\Big)_{i+\frac12,j} - v^{h,1}_{L_{i,j}} = \sum_{\Vert (r,s)\Vert \leq R} \hat{v}^1_L(r,s) e^{i\pi \frac{i r + j s}{N}} \sum_{m,n=-\nu}^\nu \eta^1_{m,n} \Big ( e^{i\pi \frac{m r + n s }{N}} - 1 \Big ) \\
    =&\sum_{\Vert (r,s)\Vert \leq R} \frac{\sum_{m,n=-\nu}^\nu \eta^1_{m,n} \big (e^{i\pi h (m r + n s)} - 1\big )}{1 + \sum_{m,n=-\nu}^\nu \eta^1_{m,n} \big (e^{i\pi h (m r + n s)} - 1 \big )}\Big (\hat{u}(r,s) \frac{e^{i\pi h r} - 1}{h} - \hat{w}^1(r,s) \Big ) e^{i\pi h(i r + j s)}.
\end{align*}
Using trigonometric identities and the triangle inequality, this can be bounded with a suitable constant $C_3$ as follows: 
\begin{align*}
    &\Big|\Big (\frac{u^h_{i+1,j}-u^h_{i,j}}{h} - w^{h,1}\Big )_{i+\frac12,j} - v^{h,1}_{L_{i,j}}\Big|  \leq  \\  & \sum_{\Vert (r,s)\Vert \leq R} \frac{\sum_{m,n=-\nu}^\nu |\eta^1_{m,n}| |\sin(\pi \frac{m r + n s}{2}h)|}{1 - \sum_{m,n=-\nu}^\nu |\eta^1_{m,n}| |\sin(\pi \frac{m r + n s}{2}h)|}\Big ( |\hat{u}(r,s)| \frac{2|\sin( \frac{\pi r h}{2})|}{h} + |\hat{w}^1(r,s)| \Big ) \leq \frac{C_3 h}{\sqrt{2}}. \notag
\end{align*}
For the second component, the same procedure leads to
\begin{align*}
     \Big|\Big (\frac{u^h_{i,j+1}-u^h_{i,j}}{h} - w^{h,2}\Big)_{i,j+\frac12} - v^{h,2}_{L_{i,j}}\Big | \leq  \frac{C_3 h}{\sqrt{2}}.
\end{align*}
Now, proceeding similarly for the second part 
in~\eqref{eq:tgv_limsup_aim} related to $v_K^h$, we want to find $v_K^{h,1}$, $v_K^{h,2}$, $v_K^{h,3}$ such that
\begin{align} \label{eq:limsup_vK}
    \bigg(\frac{w^{h,1}_{i+\frac32,j} -w^{h,1}_{i+\frac12,j}}{h} \bigg)_{i+1,j} &= \sum_{m,n=-\nu}^\nu \xi^1_{m,n}v^{h,1}_{K_{i+m,j+n}},\\    \bigg(\frac{w^{h,1}_{i+\frac12,j+1} -w^{h,1}_{i+\frac12,j}}{h} + \frac{w^{h,2}_{i+1,j+\frac12} -w^{h,2}_{i,j+\frac12}}{h} \bigg )_{i+\frac12,j+\frac12} &= \sum_{m,n=-\nu}^\nu \xi^2_{m,n}v^{h,2}_{K_{i+m,j+n}}, \label{eq:limsup_vK2} 
    \\
    \bigg(\frac{w^{h,2}_{i,j+\frac32} -w^{h,2}_{i,j+\frac12}}{h} \bigg)_{i,j+1} &= \sum_{m,n=-\nu}^\nu \xi^3_{m,n}v^{h,3}_{K_{i+m,j+n}}. \label{eq:limsup_vK3}
\end{align}
We consider $v^{h,1}_K$ as an example, for which the same procedure as for $v_L^{h,1}$ can be performed, where the discrete Fourier transform is applied to the first line in~\eqref{eq:limsup_vK}, and the filters are inverted to obtain
\begin{equation*}
\hat{v}^1_K(r,s) =   \frac{\hat{w}^1(r,s) \frac{e^{i\pi r h} -1}{h}}{1+ \sum_{m,n=-\nu}^\nu \xi^1_{m,n}\big (e^{i\pi h(m r + n s)} - 1 \big )}.
\end{equation*}
Then, $v^{h,1}_{K_{i,j}}$ is computed by using the inverse discrete Fourier transform. 
Finally, this can be used to obtain 
\begin{align*}
    & \bigg(\frac{w^{h,1}_{i+\frac32,j} -w^{h,1}_{i+\frac12,j}}{h} \bigg)_{i+1,j}- v_{K_{i,j}}^{h,1} \\
    =&\sum_{\Vert (r,s)\Vert \leq R} \frac{\sum_{m,n=-\nu}^\nu \xi^1_{m,n} \Big (e^{i\pi h(m r + n s )} - 1 \Big)}{1 + \sum_{m,n=-\nu}^\nu \xi^1_{m,n} \Big (e^{i\pi h(m r + n s )} - 1\Big )} \hat{w}^1(r,s) \frac{e^{i\pi r h } - 1}{h} e^{i\pi h (i r  + j s)} \notag,
\end{align*}
which can be bounded such that
\begin{equation*}
    \bigg |\bigg(\frac{w^{h,1}_{i+\frac32,j} -w^{h,1}_{i+\frac12,j}}{h} \bigg)_{i+1,j} - v_{K_{i,j}}^{h,1} \bigg | \leq \frac{C_4 h}{2\sqrt{2}}. 
\end{equation*}
Likewise, an analogous computation for $v_K^{h,2}$,$v_K^{h,3}$ in~\eqref{eq:limsup_vK2} and \eqref{eq:limsup_vK3} results in
\begin{align*}
    \hat{v}^2_K(r,s) =   \frac{\hat{w}^1(r,s) \frac{e^{i\pi s h} -1}{h} + \hat{w}^2(r,s) \frac{e^{i\pi r h} -1}{h}}{1+ \sum_{m,n=-\nu}^\nu \xi^2_{m,n}\big (e^{i\pi h(m r + n s)} - 1 \big )}, \\
    \hat{v}^3_K(r,s) =   \frac{\hat{w}^2(r,s) \frac{e^{i\pi s h} -1}{h}}{1+ \sum_{m,n=-\nu}^\nu \xi^3_{m,n}\big (e^{i\pi h(m r + n s)} - 1 \big )}.
\end{align*}
It follows that~\eqref{eq:limsup_general} holds with an error of $\mathcal{O}(h)$ and thus we can conclude that
\begin{align*}
 &
 \lim \sup_{h\to 0}h^2 \sum_{i,j=1}^N \alpha_1 \sqrt{\Big (\tfrac{1}{h}(u^h_{i+1,j} - u^h_{i,j})-  w^{h,1}_{i+\frac12,j}\Big)^2 + \Big (\tfrac{1}{h}(u^h_{i,j+1} - u^h_{i,j}) -  w^{h,2}_{i,j+\frac12}\Big )^2} + \\
\textstyle
& h^2  \sum_{i,j=1}^N \alpha_0  \sqrt{\tfrac{1}{h^2}( w^{h,1}_{i+\frac32,j}-w^{h,1}_{i+\frac12,j})^2 +\tfrac{1}{h^2}(w^{h,2}_{i,j+\frac32}-w^{h,2}_{i,j+\frac12})^2 \rule{0pt}{10pt}\ } \notag \\
\textstyle
& \overline{\rule{0pt}{17pt}
    + \tfrac{2}{h^2} \Big (\tfrac 1 2 (w^{h,1}_{i+\frac12,j+1} - w^{h,1}_{i+\frac12,j}) + \tfrac 1 2 (w^{h,2}_{i+1,j+\frac12} - w^{h,2}_{i,j+\frac12})\Big )^2
\ }
 \notag \\ 
\leq&  \int_\Omega \alpha_1 |\D u -w| + \alpha_0 |\Epsop w| \d{x},
\end{align*}
which readily implies~\eqref{eq:limsup_general}.

So far, we assumed that $u,w$ have a finite number of modes, which also translates to the previously defined sequence $(u_k,w_k)$.
Thanks to the smoothness, we know that $u_{k,h}\coloneqq(u_k(ih,jh))_h$, viewed as a piecewise constant function over pixels of size $h\times h$, converges (uniformly) to $u_k$. 
This means that we can again use a diagonal argument and let $h\to 0$ and $k\to\infty$ simultaneously, combining~\eqref{eq:limsup_general} and~\eqref{eq:approxfinitemodes}, to define $u^h=u_{k_h,h}$ which converges to $u$ and satisfies
\begin{equation}\label{eq:Glimsup}
    \limsup_{h\to 0} \text{TGV}^2_{\alpha,h}(u^h)
\le \text{TGV}^2_\alpha(u),
\end{equation}
showing the $\Gamma-\limsup$ property. 
This can be seen as follows: we let $h_0=1$ and for $k\ge 1$, we recursively find $h_k< h_{k-1}$ such that -- thanks to~\eqref{eq:limsup_general} -- we obtain
\[
\|u_{h,k}-u_k\|_1 \le \frac{1}{k},\qquad
\text{TGV}^2_{\alpha,h}(u_{h,k})
\le \int_\Omega\alpha_1|Du_k-w_k|+\alpha_0|E w_k| + \frac{1}{k},
\]
for $h\le h_k$. 
Then, we let $k_h = k$ whenever $h_{k+1}<h\le h_k$. 
We easily conclude that $u^h\to u$, and, thanks to~\eqref{eq:approxfinitemodes}, we get~\eqref{eq:Glimsup}.
\end{proof}

\bibliographystyle{unsrtnat}
\bibliography{references}  






\end{document}